\documentclass[a4paper,12pt]{article}

\usepackage{amsmath,amsthm,amssymb}
\allowdisplaybreaks

\usepackage{mathrsfs}

\usepackage{bbold}

\makeatletter
\DeclareFontFamily{U}{escorr@fam}{}
\DeclareFontShape{U}{escorr@fam}{m}{n}{%
 <5>sfixed * [3.75] cmsy10%
 <6>sfixed * [4.5] cmsy10%
 <7>sfixed * [5.25] cmsy10%
 <8>sfixed * [6] cmsy10%
 <9>sfixed * [6.75] cmsy10%
 <10>sfixed * [7.5] cmsy10%
 <10.95>sfixed * [8.2125] cmsy10%
 <12>sfixed * [9] cmsy10%
 <14.4>sfixed * [10.8] cmsy10%
 <17.28>sfixed * [12.96] cmsy10%
 <20.74>sfixed * [15.555] cmsy10%
 <24.88>sfixed * [18.66] cmsy10%
 }{}
\DeclareMathAlphabet{\escorr@alph}{U}{escorr@fam}{m}{n}
\newcommand*{\escorr@macro}[2]{\savebox{0}{$#1\escorr@alph{\mathchar"710D}\n@space$}
                               \savebox{1}[\wd0]{\hss$#1\mathchar"002F\n@space$\hss}
                               \dimen0=\ht1 \advance\dimen0by\dp1 \divide\dimen0by16
                               \mbox{\rlap{\usebox{1}}\raisebox{\dimen0}{\usebox{0}}}}
\DeclareRobustCommand{\Memptyset}{{\mathpalette\escorr@macro{}}}
\DeclareRobustCommand{\Mbfemptyset}{{\pmb{\pmb{\Memptyset}}}}
\makeatother

\makeatletter
\newcommand*{\bigcorr@macro}[2]{\sbox{0}{\mbox{$#1($}}\dimen0=\ht0
				\advance\dimen0 by \dp0
				\multiply\dimen0 by #2 \divide\dimen0 by 100}
\newcommand*{\bigcorr@big}[2]{\mbox{$#1\left#2\bigcorr@macro{#1}{85}\vrule
			       height \dimen0 depth 0pt width 0pt\right.\n@space$}}
\newcommand*{\bigcorr@Big}[2]{\mbox{$#1\left#2\bigcorr@macro{#1}{115}\vrule
			       height \dimen0 depth 0pt width 0pt\right.\n@space$}}
\newcommand*{\bigcorr@bigg}[2]{\mbox{$#1\left#2\bigcorr@macro{#1}{145}\vrule
			       height \dimen0 depth 0pt width 0pt\right.\n@space$}}
\newcommand*{\bigcorr@Bigg}[2]{\mbox{$#1\left#2\bigcorr@macro{#1}{175}\vrule
			       height \dimen0 depth 0pt width 0pt\right.\n@space$}}
\DeclareRobustCommand*{\big}[1]{{\mathpalette\bigcorr@big{#1}}}
\DeclareRobustCommand*{\Big}[1]{{\mathpalette\bigcorr@Big{#1}}}
\DeclareRobustCommand*{\bigg}[1]{{\mathpalette\bigcorr@bigg{#1}}}
\DeclareRobustCommand*{\Bigg}[1]{{\mathpalette\bigcorr@Bigg{#1}}}
\DeclareRobustCommand*{\bi}[1]{{#1}}
\DeclareRobustCommand*{\bil}[1]{\mathopen{\bi{#1}}}

\DeclareRobustCommand*{\bir}[1]{\mathclose{\bi{#1}}}
\makeatother

\renewcommand{\leq}{\leqslant}

\renewcommand{\geq}{\geqslant}

\newcommand{\pointir}{\unskip.\kern 0.32em---\kern 0.32em\ignorespaces}
\newtheoremstyle{mynumtheo}   
   {5pt}                      
   {5pt}                      
   {\slshape}                 
   {}                         
   {\bfseries}                
   {{\bfseries\pointir}}      
   {0pt}                      
   {$\pmb{\mbox{{\sc #1}}}$\ #2{\thmnote{\mdseries\ [#3]}}} 
\theoremstyle{mynumtheo}
\newtheorem{criterion}{Criterion}
\newtheorem{proposition}{Proposition}

\newtheoremstyle{myunnumtheo} 
   {5pt}                      
   {5pt}                      
   {\slshape}                 
   {}                         
   {\bfseries}                
   {{\bfseries\pointir}}      
   {0pt}                      
   {$\pmb{\mbox{{\sc #1}}}$\thmnote{\mdseries\ [#3]}} 
\theoremstyle{myunnumtheo}
\newtheorem{lemma}{Lemma}
\newtheorem{definition}{Definition}
\newtheorem{T}{Theorem}

\DeclareMathOperator{\parts}{\mathrm{part}}
\DeclareMathOperator{\pb}{\mathbf{P}}
\DeclareMathOperator{\qb}{\mathbf{Q}}
\DeclareMathOperator*{\Mlim}{lim\vphantom{p}}
\newcommand{\ZZ}{{\mathbb{Z}}}
\newcommand{\RR}{{\mathbb{R}}}
\newcommand{\supp}{{\mathfrak S}}
\newcommand{\wt}{\widetilde}
\newcommand{\sm}{\setminus}
\newcommand{\st}{\,:\,}
\newcommand{\cmp}{^{\mathrm{c}}}
\newcommand{\zl}{\vphantom{\wt\Lambda}}
\newcommand{\zt}{\vphantom{t}}
\newcommand*{\bm}[1]{\boldsymbol{#1}}
\newcommand*{\cm}[1]{\mathscr{#1}}
\newcommand*{\bca}[1]{\boldsymbol{\mathcal{#1}}}
\newcommand*{\ov}[1]{\overline{\bm #1}}
\newcommand*{\fov}[1]{\wt{\bm #1}}

\begin{document}

\title{On Gibbsianness of Random Fields}
\author{S.~Dachian
	\thanks{Laboratoire de Math\'ematiques, Universit\'e
		Blaise Pascal (Clermont-Ferrand 2), 63177 Aubi\`ere CEDEX,
	        FRANCE, e-mail:
                \texttt{Serguei.Dachian@math.univ-bpclermont.fr}}
	\and
	B.S.~Nahapetian
	\thanks{Institute of Mathematics, National Academy of Sciences of
		Armenia, 24-b Av.~M.~Bagramian, 375019 Yerevan, ARMENIA,
		e-mail: \texttt{nahapet@instmath.sci.am}}}
\date{}
\maketitle

\vglue 2.5cm
\begin{abstract}
\noindent
The problem of characterization of Gibbs random fields is considered. Various
Gibbsianness criteria are obtained using the earlier developed one-point
framework which in particular allows to describe random fields by means of
either one-point conditional or one-point finite-conditional
distributions. The main outcome are the criteria in terms of one-point
finite-conditional distribution, one of which can be taken as a purely
probabilistic definition of Gibbs random field.

\vglue 2cm
\centerline{\textbf{R\'esum\'e}}

\bigskip
\noindent
Le probl\`eme consid\'er\'e est celui de la caract\'erisation des champs
al\'eatoires de Gibbs. Divers crit\`eres de gibbsianit\'e sont obtenus en
utilisant l'approche uni\-ponctuelle d\'evelopp\'ee ant\'erieurement et qui
permet en particulier de d\'ecrire les champs al\'eatoires soit par des
probabilit\'es conditionnelles uniponctuelles, soit par des probabilit\'es
fini-conditionnelles uniponctuelles. Les r\'esultats principaux sont les
crit\`eres exprim\'es en termes de probabilit\'es fini-conditionnelles
uniponctuelles, l'un desquels peut \^etre pris comme une d\'efinition purement
probabiliste du champ al\'eatoire de Gibbs.

\vglue 2.5cm
\noindent
\end{abstract}
\newpage

\section*{Introduction}

The classes of processes considered in the random processes theory are usually
specified by some properties of their finite-dimensional or conditional
distributions. However in practice, the study of a particular class usually
goes through some representation theorem expressing processes in terms of
simple and convenient objects, such as transition matrices for Markov chains,
characteristic functions for processes with independent increments, spectral
functions for stationary processes, and so on.

The situation is quite different for the class of Gibbs random
fields. Historically, instead of being characterized by some properties of
their finite-dimensional or conditional distributions, Gibbs random fields
have been defined directly by representation of their conditional
distributions in terms of potentials. And only afterwards the problem of
probabilistic characterization of Gibbs random fields was considered.

It was shown in Kozlov~\cite{kozl} and Sullivan~\cite{sull2} that Gibbs random
fields (with uniformly convergent potentials) can be characterized by strict
positivity and quasilocality of their conditional distributions. More
precisely, in order for a random field to be Gibbsian, its conditional
distribution (which consists of conditional probabilities on finite volumes
with conditions on the entire exterior and is defined up to a set of
probability zero) must have a version which is a strictly positive quasilocal
specification. As we see, this criterion imposes conditions on an object
(conditional distribution) which is neither unambiguously defined, nor
constructive (is infinite-dimensional), and, moreover, does not determine the
random field uniquely (phase transitions). In our opinion, it is preferable
that a characterization be in terms of an object which does not have these
features.

As a matter of fact, such characterization already exists for the particular
class of real-valued finite-range potentials. It was shown in
Averintsev~\cite{aver,aver2,aver3} and Sullivan~\cite{sull} that Gibbs random
fields (with such potentials) can be characterized by strict positivity and
Markov properties. Note that for strictly positive random fields the Markov
property can be formulated using only conditional probabilities on single
sites with finite-volume conditions $\bigl($see, for example,
Suomela~\cite{suom}$\bigr)$. These probabilities are defined unambiguously and
in constructive manner as ratios of finite-dimensional
probabilities. Moreover, according to Dalalyan and Nahapetian~\cite{DaN}, the
system of all such probabilities, which we call one-point finite-conditional
distribution, uniquely determines (can be identified with) the random field.

The aim of this paper is to characterize Gibbs random fields by some
properties of their one-point finite-conditional distributions in the general
case of uniformly convergent potentials. It is worth mentioning that such
characterization is very natural in light of and was made possible due to the
one-point framework developed in some recent papers. Namely, an approach
towards description of random fields was developed by the authors
in~\cite{DN2,DN3,DN4}, where Dobrushin's well-known description of random
fields by means of conditional distributions was reduced to description of
random fields by means of one-point conditional distributions (the system of
conditional probabilities on single sites with conditions on the entire
exterior). Later on, a closely related and in some way complementary
description of random fields based on one-point finite-conditional
distributions was proposed in Dalalyan and Nahapetian~\cite{DaN}.

The main outcome of the present work are random field Gibbsianness criteria in
terms of one-point finite-conditional distribution. Let us emphasize that one
of these criteria can be taken as a purely probabilistic definition of Gibbs
random field, turning the usual definition into a representation theorem for
Gibbs random field's conditional distribution. The criteria are established in
Section~\ref{GC}, which also contains some additional results on
characterization of Gibbsian specifications and Gibbs random fields, along
with a brief survey of the known ones. Before that, some necessary notation
and prerequisites are given in Section~\ref{P}, while the above mentioned
one-point framework is presented in Section~\ref{OPF}.

\section{Preliminaries}
\label{P}

In this section we briefly recall some necessary notions and facts from the
theory of Gibbs random fields.

\subsection{Random fields}

We consider random fields on the $\nu$-dimensional integer lattice $\ZZ^\nu$,
that is, probability measures $\pb$ on $\bigl(\cm X^{\ZZ^\nu},\cm
F^{\ZZ^\nu}\bigr)$ where $(\cm X,\cm F)$ is some measurable space of values on
single sites (\emph{state space\/}). Usually the space $\cm X$ is assumed to
be endowed with some topology~$\cm T$, and $\cm F$ is assumed to be the Borel
$\sigma$-algebra for this topology. In this work we concentrate on the case
when $\cm X$ is finite,~$\cm T$ is the discrete topology, and~$\cm F$ is the
total $\sigma$-algebra, that is, $\cm F=\cm T=\parts(\cm X)$.

For any $S\subset\ZZ^\nu$, we denote $\cm E(S)$ the set of all finite subsets
of $S$, that is, we put $\cm E(S)=\bigl\{\Lambda\subset
S\st\bil|\Lambda\bir|<\infty\bigr\}$ where $\bil|\Lambda\bir|$ is the number
of points of the set $\Lambda$. For convenience of notation we will omit
braces for one-point sets, that is, will write~$t$ instead of $\{t\}$. We put
also $\cm E^*(S)=\cm E\sm\{\Memptyset\}$. For $S=\ZZ^\nu$ we write $\cm E=\cm
E(\ZZ^\nu)$ and~$\cm E^*=\cm E^*(\ZZ^\nu)$.

For any $S\subset\ZZ^\nu$, the space ${\cm X}^S$ is the space of all
configurations on $S$. If $S=\Memptyset$, we assume that the space $\cm
X^\Memptyset=\{\Mbfemptyset\}$ where $\Mbfemptyset$ is the empty
configuration. For any $T,S\subset\ZZ^\nu$ such that $T\subset S$ and any
configuration $\bm x=\{x_t,\ t\in S\}$ on $S$, we denote $\bm x_{\zl T}$ the
\emph{subconfiguration\/} (\emph{restriction\/}) of $\bm x$ on~$T$ defined by
$\bm x_{\zl T}=\{x_t,\ t\in T\}$. For any $T,S\subset\ZZ^\nu$ such that $T\cap
S=\Memptyset$ and any configurations $\bm x$ on $T$ and $\bm y$ on $S$, we
denote $\bm x\bm y$ the \emph{concatenation\/} of $\bm x$ and $\bm y$, that
is, the configuration on $T\cup S$ equal to $\bm x$ on $T$ and to $\bm y$ on
$S$. For any configuration $\bm x\in{\cm X}^S$, the set $S\subset\ZZ^\nu$ will
be called \emph{support\/} of $\bm x$ and we will write $S=\supp(\bm x)$. For
any $\Lambda\in\cm E$, we denote
$$
\wt{{\cm X}^\Lambda}=\bigcup_{\wt\Lambda\in\cm E^*(\Lambda\cmp
)}{\cm X}^{\wt\Lambda}
$$
the space of all configurations with non-empty finite support contained in the
exterior of $\Lambda$.

For any $S\subset\ZZ^\nu$, a probability distribution on $\cm X^S$ will be
denoted by $\pb_S$. Note that if $S=\Memptyset$ there exists only one
probability distribution $\pb_\Memptyset(\Mbfemptyset)=1$.  For any
$T,S\subset\ZZ^\nu$ such that $T\subset S$ and any $\pb_S$, we denote
${\bigl(\pb_S\bigr)}_T$ the \emph{marginal distribution\/}
(\emph{restriction\/}) of $\pb_S$ on $T$. If $\Lambda\in\cm E$ and
$I\subset\Lambda$, we can write $\pb_\Lambda=\bigl\{\pb_\Lambda(\bm x),\ \bm
x\in\cm X^\Lambda\bigr\}$ and
$$
{\bigl(\pb_\Lambda\bigr)}_I (\bm x)=\sum_{\bm y\in\cm X^{\Lambda\sm I}}
\pb_\Lambda(\bm x\bm y),\quad \bm x\in\cm X^I.
$$

Any random field $\pb$ on $\ZZ^\nu$ is uniquely determined by (can be
identified with) the system $\{\pb_\Lambda,\ \Lambda\in\cm E\}$ of its
\emph{finite-dimensional distributions\/} which are consistent in the sense
that for any\/ $\Lambda\in\cm E$ and\/ $I\subset\Lambda$ we have\/
${\bigl(\pb_\Lambda\bigr)}_I=\pb_I$.

Finally, a random field $\pb$ will be called \emph{strictly positive\/} if for
any $\Lambda\in\cm E$ the finite-dimensional distribution $\pb_\Lambda$ is
\emph{strictly positive}, that is, $\pb_\Lambda(\bm x)>0$ for all~$\bm x\in\cm
X^\Lambda$.

\subsection{Finite-conditional and conditional distributions\\
            of random fields}

Let $\pb$ be some random field. For any $\Lambda\in\cm E$, we denote $\pb_{\cm
E^*(\Lambda\cmp)}$ the measure on $\wt{\cm X^\Lambda}$ whose projection on
${\cm X}^{\wt\Lambda}$ is $\pb_{\wt\Lambda}$ for any $\wt\Lambda\in\cm
E^*(\Lambda\cmp)$, that is, $\pb_{\cm E^*(\Lambda\cmp)}$ is the direct sum of
the measures $\pb_{\wt\Lambda}$.

For all $\Lambda\in\cm E$, the ratios
$$
\bm q_\Lambda^{\fov x}(\bm x)=\frac{\pb_{\Lambda\cup\supp(\fov x)}(\bm x\fov
x)}{\pb_{\supp(\fov x)}(\fov x)}\ ,\quad \bm x\in\cm X^\Lambda,
$$
exist for $\pb_{\cm E^*(\Lambda\cmp)}$-almost all $\fov x\in\wt{{\cm
X}^\Lambda}$. Any system
$$
\wt{\bca Q}=\Bigl\{\qb_\Lambda^{\fov x},\quad \Lambda\in\cm E\text{ and }\fov
x\in\wt{\cm X^\Lambda}\Bigr\}
$$
of probability distributions such that for every $\Lambda\in\cm E$ we have
$\qb_\Lambda^{\fov x}=\bm q_\Lambda^{\fov x}$ for $\pb_{\cm
E^*(\Lambda\cmp)}$\nobreakdash-almost all $\fov x\in\wt{\cm X^\Lambda}$ will
be called \emph{finite-conditional distribution\/} of the random
field~$\pb$. The subsystem of $\wt{\bca Q}$ consisting of single-site
distributions $\bigl(\left|\Lambda\right|=1\bigr)$ will be called
\emph{one-point finite-conditional distribution\/} of~$\pb$. Note that in
general a random field may have many \emph{versions\/} both of
finite-conditional and one-point finite-conditional distributions. However,
for strictly positive random fields these distributions are uniquely
determined and consist of strictly positive elements. If this positivity is
uniform with respect to $\fov x$, the (one-point) finite-conditional
distribution will be called \emph{uniformly nonnull\/}.

Further, for all $\Lambda\in\cm E$, the limits
$$
\bm q_\Lambda^{\ov x}(\bm x)=\Mlim_{\wt\Lambda\uparrow\ZZ^\nu\sm\Lambda}
\bm q_\Lambda^{\ov x_{\wt\Lambda}}(\bm x)\ ,\quad \bm x\in\cm X^\Lambda,
$$
exist for $\pb_{\Lambda\cmp}$-almost all $\ov x\in\cm X^{\Lambda\cmp}$. Any
system
$$
\bca Q=\Bigl\{\qb_\Lambda^{\ov x},\quad \Lambda\in\cm E\text{ and }\ov x\in\cm
X^{\Lambda\cmp}\Bigr\}
$$
of probability distributions such that for every $\Lambda\in\cm E$ we have
$\qb_\Lambda^{\ov x}=\bm q_\Lambda^{\ov x}$ for
$\pb_{\Lambda\cmp}$\nobreakdash-almost all $\ov x\in\cm X^{\Lambda\cmp}$ will
be called \emph{conditional distribution\/} of the random field~$\pb$. The
subsystem of $\bca Q$ consisting of single-site distributions will be called
\emph{one-point conditional distribution\/} of~$\pb$. Note that in general a
random field $\pb$ may have many \emph{versions\/} both of conditional and
one-point conditional distributions (even if $\pb$ is strictly positive). Note
also that if a random field $\pb$ has a strictly positive version of
conditional distribution (a version all of whose elements are strictly
positive), then $\pb$ is necessarily strictly positive itself.

Concluding this section let us remark that (one-point) conditional
distribution can be deduced from (one-point) finite-conditional distribution,
but not the other way around. So, the latter contains ``more complete''
information about the random field than the former. Especially it becomes
apparent in the Markov case, when (one-point) conditional distribution can be
considered as a subsystem of (one-point) finite-conditional
distribution. Indeed, let $\pb$ be a Markov random field and let
$\partial\Lambda$ denote the neighborhood of the set $\Lambda$. As we have
$\qb_\Lambda^{\ov x_{\zl}}=\qb_\Lambda^{\ov x_{\zl\partial\Lambda}}$, the
elements of (one-point) conditional distribution of $\pb$ can also be
considered as elements of (one-point) finite-conditional distribution of
$\pb$. The converse is not true since not all the elements of the latter
correspond to the elements of the former, but only the elements
$\qb_\Lambda^{\fov x}$ such that $\supp(\fov x)\supset\partial\Lambda$.

\subsection{Description of random fields\\
            by means of conditional distributions}

The well-known description of random fields by means of conditional
distributions introduced by Dobrushin in~\cite{dobr,dobr2,dobr3} is carried
out in terms of specifications. A system
$$
\bca Q=\Bigl\{\qb_\Lambda^{\ov x},\quad \Lambda\in\cm E\text{ and }\ov x\in\cm
X^{\Lambda\cmp}\Bigr\}
$$
of probability distributions is called \emph{specification\/} if
\begin{equation}
\label{spec}
\begin{gathered}
\qb_{\Lambda}^{\ov x}(\bm x\bm y)={\bigl(\qb^{\ov
x}_{\Lambda}\bigr)}_{\Lambda\sm I}(\bm x)\,\qb_I^{\ov x\bm x}(\bm y)\\
\text{for all~~$\Lambda\in\cm E$,~~$I\subset\Lambda$,~~$\bm x\in\cm
X^{\Lambda\sm I}$,~~$\bm y\in\cm X^I$~~and~~$\ov x\in\cm X^{\Lambda\cmp}$.}
\end{gathered}
\end{equation}
Note that any version of conditional distribution of a random field $\pb$
satisfies a somewhat weaker than~\eqref{spec} condition, where
$\pb_{\Lambda\cmp}$-almost all (and not necessarily all) $\ov x\in\cm
X^{\Lambda\cmp}$ are considered. However, any random field possesses at least
one version of conditional distribution being a specification $\bigl($see
Goldstein~\cite{gold}, Preston~\cite{pres} and Sokal~\cite{soka}$\bigr)$.

One of the main goals of Dobrushin's theory is to study the set of all random
fields \emph{compatible\/} with a given specification, that is, having it as a
version of conditional distribution. The best-known sufficient conditions for
existence and for uniqueness of random fields compatible with a given
specification are quasilocality and Dobrushin's uniqueness conditions
respectively. The first one will play an important role in our considerations,
so we recall it below.

Let $S\subset\ZZ^\nu$. A real-valued function $g$ on $\cm X^S$ is called
\emph{quasilocal\/} if
$$
\Mlim_{\Lambda\uparrow\vphantom{\cm X^S}S}\ \ \sup_{\bm x,\bm y\in\cm
X^S\vphantom{\cm X^S}\st\bm x_{\zl\Lambda}=\bm y_{\zl\Lambda}}\ \ \bigl|g(\bm
x)-g(\bm y)\bigr|=0,
$$
or equivalently if $g$ is a uniform limit of functions depending only on
values of configuration on finite sets of sites (\emph{local\/}
functions). Note also that the quasilocality is nothing but continuity with
respect to the topology $\cm T^S$ and, taking into account that $\cm X^S$ is
compact, the strict positivity and uniform nonnullness conditions are
equivalent for quasilocal functions.

A specification $\bca Q=\bigl\{\qb_\Lambda^{\ov x},\quad \Lambda\in\cm E\text{
and }\ov x\in\cm X^{\Lambda\cmp}\bigr\}$ is called \emph{(quasi)local\/} if
for any $\Lambda\in\cm E$ and $\bm x\in\cm X^\Lambda$ the function $\ov
x\mapsto\qb_\Lambda^{\ov x}(\bm x)$ on~$\cm X^{\Lambda\cmp}$ is (quasi)local.

Finally, a specification will be called \emph{strictly positive\/} if all its
elements are strictly positive.

\subsection{Gibbs random fields and Gibbsian specifications}

The main object of consideration of the present paper are Gibbs random
fields. The latters are defined in terms of Gibbsian specifications, which in
turn are defined in terms of potentials.

Any function $\Phi$ on $\wt{\cm X^\Memptyset}$ taking values in
$\RR\cup\{+\infty\}$ is called \emph{(interaction) potential\/}. A potential
$\Phi$ is called \emph{convergent\/} if it is real-valued and the series
\begin{equation}
\label{convpotlims}
\sum_{\wt J\in\cm E(t\cmp)}\Phi\bigl(x\ov x_{\wt J}\bigr)
\end{equation}
converge for all $t\in\ZZ^\nu$, $x\in\cm X^t$ and $\ov x\in\cm
X^{t\cmp}$.

A potential $\Phi$ is called \emph{uniformly convergent\/} if it is convergent
and the convergence in~\eqref{convpotlims} is uniform with respect to $\ov x$.

A potential $\Phi$ is called \emph{finite-range potential\/} if for any
$t\in\ZZ^\nu$ there exist only a finite number of sets $\wt J\in\cm E(t\cmp)$
such that $\Phi\not\equiv 0$ on $\cm X^{t\cup\wt J}$. Note that any
real-valued finite-range potential is uniformly convergent.

For an arbitrary convergent potential $\Phi$ one can construct the
specification \hbox{$\bca Q=\bigl\{\qb_\Lambda^{\ov x},\quad \Lambda\in\cm
E\text{ and }\ov x\in\cm X^{\Lambda\cmp}\bigr\}$} given by Gibbs formulae
\begin{equation}
\label{gibbs}
\qb_\Lambda^{\ov x}(\bm x)=\frac {\exp\bigl(-U_\Lambda^{\ov x}(\bm x)\bigr)}
{\sum\limits_{\bm y\in\cm X^\Lambda} \exp\bigl(-U_\Lambda^{\ov x}(\bm
y)\bigr)}\ ,\quad\Lambda\in\cm E,\ \bm x\in\cm X^\Lambda,\ \ov x\in\cm
X^{\Lambda\cmp},
\end{equation}
where
\begin{equation}
\label{energy}
U_\Lambda^{\ov x}(\bm x)\ =\!\!\sum_{\zl J\st\Memptyset\ne J\subset\Lambda}\ \
\sum_{\wt J\in\cm E(\Lambda\cmp)}\Phi\bigl(\bm x_{\zl J}\ov x_{\wt
J}\bigr),\quad \Lambda\in\cm E,\ \bm x\in\cm X^\Lambda,\ \ov x\in\cm
X^{\Lambda\cmp}.
\end{equation}
The specification $\bca Q$ is called \emph{Gibbsian with potential\/} $\Phi$.
Any random field compatible with $\bca Q$ is called \emph{Gibbs random field
with potential\/} $\Phi$.

In this paper we consider uniformly convergent potentials only, so Gibbsian
specifications and Gibbs random fields with uniformly convergent potentials
will be called shortly \emph{Gibbsian specifications\/} and \emph{Gibbs random
fields\/} correspondingly.

\section{One-point framework}
\label{OPF}

The idea that it is possible to describe and study random fields using only
one-point conditional probabilities goes back to Dobrushin~\cite{dobr}. It was
realized in the authors' works~\cite{DN2,DN3,DN4} using one-point conditional
distributions, and in Dalalyan and Nahapetian~\cite{DaN} using one-point
finite-conditional distributions.

\subsection{Description of specifications and random fields\\
            by means of one-point conditional distributions}
\label{OneSpec}

In this section we briefly recall the main results of the authors'
works~\cite{DN2,DN3,DN4}.

In these papers, under wide positivity assumptions (\emph{very weak
positivity\/}) a necessary and sufficient conditions for a system
$\bigl\{\qb_t^{\ov x},\quad t\in\ZZ^\nu\text{ and }\ov x\in\cm
X^{t\cmp}\bigr\}$ of probability distributions to be contained in some
specification were established. A system satisfying these conditions was
called 1-specification. It was equally shown that the specification containing
the given 1-specification is uniquely determined by some explicit formulae
involving only the elements of this 1-specification. Moreover, since these
formulae make use of finite number of elementary operations, the entire
specification is quasilocal if and only if the 1-specification is, and the set
of random fields compatible with the 1-specification coincides with the set of
random fields compatible with the entire specification. So, whole Dobrushin's
theory can be reformulated in terms of 1-specifications and one can speak
about description of random fields by means of one-point conditional
distributions.

Let us now give some more details in the particular strictly positive case.

The definition of strictly positive 1-specification can be formulated in the
following way: a system
$$
\bca Q=\Bigl\{\qb_t^{\ov x},\quad t\in\ZZ^\nu\text{ and }\ov x\in\cm
X^{t\cmp}\Bigr\}
$$
of strictly positive probability distributions will be called
\emph{1-specification\/} if
\begin{equation}
\label{cycle}
\begin{gathered}
\qb_t^{\ov x v}(x)\,
\qb_{s\zt}^{\ov x x}(y)\,
\qb_t^{\ov x y}(u)\,
\qb_{s\zt}^{\ov x u}(v)
=
\qb_{s\zt}^{\ov x u}(y)\,
\qb_t^{\ov x y}(x)\,
\qb_{s\zt}^{\ov x x}(v)\,
\qb_t^{\ov x v}(u)\\
\text{for all~~$t,s\in\ZZ^\nu$,~~$x,u\in\cm X^t$,~~$y,v\in\cm X^s$~~and~~$\ov
x\in\cm X^{\{t,s\}\cmp}$.}
\end{gathered}
\end{equation}
Further, a 1-specification $\bca Q=\bigl\{\qb_t^{\ov x},\quad
t\in\ZZ^\nu\text{ and }\ov x\in\cm X^{t\cmp}\bigr\}$ is called
\emph{(quasi)local\/} if for any $t\in\ZZ^\nu$ and $x\in\cm X^t$ the function
$\ov x\mapsto\qb_t^{\ov x}(x)$ on~$\cm X^{t\cmp}$ is (quasi)local. Finally, a
random fields $\pb$ is called \emph{compatible\/} with a 1-specification if
the latter is a version of one-point conditional distribution of~$\pb$.

The above mentioned explicit formulae determining the elements of the
specification $\bca Q=\bigl\{\qb_\Lambda^{\ov x},\quad \Lambda\in\cm E\text{
and }\ov x\in\cm X^{\Lambda\cmp}\bigr\}$ containing the given strictly
positive 1-specification have the following form: for all $\Lambda\in\cm E$
and $\ov x\in\cm X^{\Lambda}$ one has
$$
\qb_{\Lambda}^{\ov x}(\bm x)=
\frac{\qb_{t_1}^{\ov x \bm u_{\zl\{t_2,\ldots,t_n\}}}(x_{t_1})\,
\qb_{t_2}^{\ov x \bm x_{\zl\{t_1\}} \bm u_{\zl\{t_3,\ldots,t_n\}}}(x_{t_2})
\,\cdots\,
\qb_{t_n}^{\ov x \bm x_{\zl\{t_1,\ldots,t_{n-1}\}}}(x_{t_n})}
{\qb_{t_1}^{\ov x \bm u_{\zl\{t_2,\ldots,t_n\}}}(u_{t_1})\,
\qb_{t_2}^{\ov x \bm x_{\zl\{t_1\}} \bm u_{\zl\{t_3,\ldots,t_n\}}}(u_{t_2})
\,\cdots\,
\qb_{t_n}^{\ov x \bm x_{\zl\{t_1,\ldots,t_{n-1}\}}}(u_{t_n})}\times C,
\quad\!\bm x\!\in\!\cm X^{\Lambda},
$$
where C is the normalizing factor. Here some fixed configuration $\bm u\in\cm
X^{\Lambda}$ and some enumeration $t_1,\ldots,t_n$ of elements of $\Lambda$
are chosen arbitrary. Note that the right hand side of these formulae does not
depend on this choice (correctness of the formulae) thanks to consistency
condition~\eqref{cycle}. Note also, that these formulae imply that the
specification containing a strictly positive 1-specification is necessarily
strictly positive itself.

\subsection{Description of random fields\\
            by means of one-point finite-conditional distributions}

Now we turn to the problem of description of random fields by means of
one-point finite-conditional distributions considered in Dalalyan and
Nahapetian~\cite{DaN}. This description is closely related (and in some way
complementary) to the one presented in the previous section.

First, let us note that the necessary and sufficient conditions for a system
$\wt{\bm q}=\bigl\{\qb_t^{\fov x},\quad t\in\ZZ^\nu\text{ and }\fov
x\in\wt{\cm X^t}\bigr\}$ of probability distributions to be contained in
some~system $\wt{\bca Q}=\bigl\{\qb_\Lambda^{\fov x},\quad\Lambda\in\cm
E\text{ and }\fov x\in\wt{\cm X^\Lambda}\bigr\}$ of probability distributions
satisfying
\begin{equation}
\label{fspec}
\begin{gathered}
\qb_{\Lambda}^{\fov x}(\bm x\bm y)=\qb^{\fov x}_{\Lambda\sm I}(\bm
x)\,\qb_I^{\fov x\bm x}(\bm y)\\
\text{for all~~$\Lambda\in\cm E$,~~$I\subset\Lambda$,~~$\bm x\in\cm
X^{\Lambda\sm I}$,~~$\bm y\in\cm X^I$~~and~~$\fov x\in\wt{\cm X^\Lambda}$}
\end{gathered}
\end{equation}
are the following:
\begin{equation}
\label{fcycle}
\begin{gathered}
\qb_t^{\fov x}(x)\,
\qb_{s\zt}^{\fov x x}(y)
=
\qb_{s\zt}^{\fov x}(y)\,
\qb_t^{\fov x y}(x)\\
\text{for all~~$t,s\in\ZZ^\nu$,~~$x\in\cm X^t$,~~$y\in\cm X^s$~~and~~$\fov
x\in\wt{\cm X^{\{t,s\}}}$.}
\end{gathered}
\end{equation}

Note also that if $\wt{\bm q}$ is the one-point finite-conditional
$\bigl(\wt{\bca Q}$ is the finite-conditional$\bigr)$ distribution of some
strictly positive random field, then it necessarily satisfies the
condition~\eqref{fcycle} $\bigl($the condition~\eqref{fspec}$\bigr)$. However,
in order for a strictly positive system $\wt{\bm q}$ satisfying~\eqref{fcycle}
$\bigl(\wt{\bca Q}$ satisfying~\eqref{fspec}$\bigr)$ to be the one-point
finite-conditional $\bigl($the finite-condi\-tional$\bigr)$ distribution of
some strictly positive random field one needs some additional conditions. It
turns out that such conditions are the following:
\begin{equation}
\label{fcycle2}
\begin{gathered}
\qb_t^{v}(x)\,
\qb_{s\zt}^{x}(y)\,
\qb_t^{y}(u)\,
\qb_{s\zt}^{u}(v)
=
\qb_{s\zt}^{u}(y)\,
\qb_t^{y}(x)\,
\qb_{s\zt}^{x}(v)\,
\qb_t^{v}(u)\\
\text{for all~~$t,s\in\ZZ^\nu$,~~$x,u\in\cm X^t$~~and~~$y,v\in\cm X^s$.}
\end{gathered}
\end{equation}

More precisely, in~\cite{DaN} it was shown that the strict positivity of
elements and the fulfillment of the conditions~\eqref{fcycle} and
\eqref{fcycle2} are necessary and sufficient for a system $\bigl\{\qb_t^{\fov
x},\quad t\in\ZZ^\nu\text{ and }\fov x\in\wt{\cm X^t}\bigr\}$ of probability
distributions to be the one-point finite-conditional distribution of some
strictly positive random field. It was equally shown that this random field is
uniquely determined by this system. In particular, a strictly positive random
field is uniquely determined by (can be identified with) its one-point
finite-conditional distribution, and so one can speak about description of
random fields by means of one-point finite-conditional distributions.

\section{Gibbsianness criteria}
\label{GC}

In this section we turn to the main subject of the present work: the problem
of probabilistic characterization of Gibbs random fields. But first, let us
mention several results concerning Gibbsianness of specifications. These
results are useful since taking into account the definition of Gibbs random
field, one can transform them into random field Gibbsianness criteria
in terms of conditional distribution.

\subsection{Specification Gibbsianness criteria}

The problem of characterization of the class of Gibbsian specifications with
potentials satisfying some given conditions was subject of consideration of
many authors: one can refer to Averintsev~\cite{aver,aver2,aver3} and
Sullivan~\cite{sull} for real-valued finite-range potentials,
Kozlov~\cite{kozl} and Sullivan~\cite{sull2} for uniformly convergent
potentials, our works~\cite{DN2,DN3} for more general potentials (which in
particular can assume the value~$+\infty$).

Concerning Gibbsian specifications (with uniformly convergent potentials), the
best-known criterion is the following one $\bigl($see, for example,
Georgii~\cite{geor}$\bigr)$.

\begin{criterion}[Kozlov-Sullivan]
\label{C1}
A specification is Gibbsian if and only if it is quasilocal and strictly
positive.
\end{criterion}

Combining this criterion with the results of Section~\ref{OneSpec}, one
clearly gets the following characterization already obtained by the authors
in~\cite{DN2,DN3}.

\begin{criterion}
\label{C2}
A specification is Gibbsian if and only if the 1-specification contained in it
is quasilocal and strictly positive.
\end{criterion}

Since the uniform convergence of potential assures the quasilocality of the
1-specification expressed by Gibbs formulae, one can also obtain the following
corollary of Criterion~\ref{C2}.

\begin{criterion}
\label{C3}
A specification is Gibbsian if and only if the 1-specification contained in it
can be expressed by Gibbs formulae~\eqref{gibbs} and~\eqref{energy} with some
uniformly convergent potential.
\end{criterion}

Concerning Gibbsian specifications with real-valued finite-range potentials,
let us recall that they are characterized by strict positivity and locality of
their single-site parts, and so Gibbs random fields with such potentials are
characterized by strict positivity and Markov properties. Since for strictly
positive random fields the Markov property can be formulated using only
conditional probabilities on single sites with finite-volume conditions
$\bigl($see, for example, Suomela~\cite{suom}$\bigr)$, one clearly has a
characterization of Gibbs random fields with real-valued finite-range
potentials in terms of one-point finite-conditional distribution.

Establishment of a similar characterization in the general case of uniformly
convergent potentials is not so straightforward and will be accomplished in
Sections~\ref{MainOut1} and~\ref{MainOut2}. Before that, random field
Gibbsianness criteria in terms of conditional and one-point conditional
distribution are obtained in the next two sections by means of transformation
and subsequent improvement of Criteria~\ref{C1}--\ref{C3}.

\subsection{Random field Gibbsianness criteria\\
            in terms of conditional distribution}

Combining the definition of Gibbs random field with Criterion~\ref{C1} one
gets the following well-known characterization: a random field is a Gibbs
random field if and only if it has a version of conditional distribution which
is a strictly positive quasilocal specification. This criterion can be
improved in the following way.

\begin{criterion}
\label{C4}
A random field is a Gibbs random field if and only if it has a version of
conditional distribution which is quasilocal and strictly positive.
\end{criterion}

Since the strict positivity of a version of conditional distribution implies
the strict positivity of the random field, the criterion is immediately
deduced from the following proposition which is of general interest.

\begin{proposition}
\label{P1}
If a strictly positive random field has a quasilocal version of conditional
distribution, the latter is unique and is necessarily a specification.
\end{proposition}

\begin{proof}
Let $\pb$ be a strictly positive random field. First, note that the
measure~$\pb$ is everywhere dense, that is, $\pb(A)>0$ for any non-empty open
set $A\in\cm T^{\ZZ^\nu}\sm\{\Memptyset\}$. Indeed, since such a set $A$
necessarily contains a non-empty cylinder subset $A'$, which in turn contains
a subset $\bigl\{\ov x\in\cm X^{\ZZ^\nu}\st\ov x_\Lambda=\bm x^\circ\bigr\}$
where $\Lambda\in\cm E^*$ and $\bm x^\circ\in\cm X^\Lambda$, we have
$\pb(A)\geq\pb(A')\geq\pb_\Lambda(\bm x^\circ)>0$. An important evident
property of everywhere dense measures is the following: if a continuous
function is equal to zero almost everywhere (with respect to such a measure),
then it is equal to zero everywhere.

Now, suppose $\bigl\{\qb_\Lambda^{\ov x},\quad \Lambda\in\cm E\text{ and }\ov
x\in\cm X^{\Lambda\cmp}\bigr\}$ and $\bigl\{\bm q_\Lambda^{\ov x},\quad
\Lambda\in\cm E\text{ and }\ov x\in\cm X^{\Lambda\cmp}\bigr\}$ are two
quasilocal versions of conditional distribution of $\pb$. Hence, for any
$\Lambda\in\cm E$ and $\bm x\in\cm X^\Lambda$, the function $\ov
x\mapsto\qb_\Lambda^{\ov x}(\bm x)-\bm q_\Lambda^{\ov x}(\bm x)$ on $\cm
X^{\Lambda\cmp}$ is quasilocal and equal to zero $\pb_{\Lambda\cmp}$-almost
everywhere. Since quasilocality is nothing but continuity and the measure
$\pb_{\Lambda\cmp}$ is everywhere dense, this function is equal to zero
everywhere. So, the uniqueness is proved.

Finally, suppose $\bca Q=\bigl\{\qb_\Lambda^{\ov x},\quad \Lambda\in\cm
E\text{ and }\ov x\in\cm X^{\Lambda\cmp}\bigr\}$ is (the unique) quasilocal
version of conditional distribution of $\pb$. For any $\Lambda\in\cm E$,
$I\subset\Lambda$, $\bm x\in\cm X^{\Lambda\sm I}$ and $\bm y\in\cm X^I$
consider the function $\ov x\mapsto\qb_{\Lambda}^{\ov x}(\bm x\bm
y)-{\bigl(\qb^{\ov x}_{\Lambda}\bigr)}_{\Lambda\sm I}(\bm x)\qb_I^{\ov x\bm
x}(\bm y)$ on $\cm X^{\Lambda\cmp}$. This function is clearly quasilocal and,
as it follows from the properties of conditional probabilities, is equal to
zero $\pb_{\Lambda\cmp}$-almost everywhere. Hence it is equal to zero
everywhere, and so $\bca Q$ is a specification.
\end{proof}

Let us note that Criterion~\ref{C4} was as a matter of fact obtained in
Sullivan~\cite{sull2} using a different approach.

\subsection{Random field Gibbsianness criteria\\
            in terms of one-point conditional distribution}

Criterion~\ref{C4} characterizes Gibbs random fields in terms of conditional
distribution. However, in view of Section~\ref{OneSpec}, it should be possible
to do it in terms of one-point conditional distribution.

Indeed, combining the definition of Gibbs random field with Criterion~\ref{C2}
and taking into account the results of Section~\ref{OneSpec}, one gets the
following characterization: a random field is a Gibbs random field if and only
if it has a version of one-point conditional distribution which is a strictly
positive quasilocal 1-specification. As in the preceding section we can
improve this criterion in the following way.

\begin{criterion}
\label{C5}
A random field is a Gibbs random field if and only if it has a version of
one-point conditional distribution which is quasilocal and strictly positive.
\end{criterion}

The criterion is immediately deduced from the following two propositions which
are of general interest.

\begin{proposition}
If a random field $\pb$ has a strictly positive version of one-point
conditional distribution, then $\pb$ is strictly positive itself.
\end{proposition}

\begin{proof}
Let us suppose that the random field $\pb$ is not strictly positive. In this
case we can find some $\Lambda\in\cm E^*$, $t\in\Lambda$ and $\bm z\in\cm
X^\Lambda$ such that $\pb_{\Lambda}(\bm z)=0$ and $\pb_{\Lambda\sm t}(\bm
z_{\Lambda\sm t})>0$ $\bigl($recall that
$\pb_\Memptyset(\Mbfemptyset)=1\bigr)$. Now denote
$$
A=\Bigl\{\ov x\in\cm X^{t\cmp}\st\ov x_{\zl\Lambda\sm t}=\bm z_{\zl\Lambda\sm
t}\Bigr\}\:.
$$
Obviously $\pb_{t\cmp}(A)=\pb_{\zl\Lambda\sm t}(\bm z_{\zl\Lambda\sm
t})>0$. Introduce also
$$
B=\bigcap_{\wt\Lambda\in\cm E(t\cmp)}\Bigl\{\ov x\in\cm X^{t\cmp}\st
\pb_{\wt\Lambda}(\ov x_{\wt\Lambda})>0\Bigr\}\:.
$$
Since $B$ is a countable intersection of sets of probability~$1$, we have
$\pb_{t\cmp}(B)=1$. So, it comes $\pb_{t\cmp}(A\cap B)>0$.

For all $\ov x\in A\cap B$ and all $\wt\Lambda\in\cm E(t\cmp)$ such that
$\wt\Lambda\supset\Lambda\sm t$, we have
$$
\bm q_t^{\ov x_{\wt\Lambda}}(\bm z_{\zl
t})\triangleq\frac{\pb_{t\cup\wt\Lambda}(\bm z_{\zl t}\ov x_{\wt\Lambda})}
{\pb_{\wt\Lambda}(\ov x_{\wt\Lambda})}=0.
$$
Hence, for all $\ov x\in A\cap B$ we get
$$
\Mlim_{\wt\Lambda\uparrow\ZZ^\nu\sm t}\bm q_t^{\ov x_{\wt\Lambda}}(\bm z_{\zl
t})=0
$$
which contradicts the existence of a strictly positive version of one-point
conditional distribution of $\pb$.
\end{proof}

\begin{proposition}
If a strictly positive random field has a quasilocal version of one-point
conditional distribution, the latter is unique and is necessarily a
1-specification.
\end{proposition}

\begin{proof}
The uniqueness is proved following exactly the same argument as in the proof
of Proposition~\ref{P1}.

To prove the second assertion, suppose $\bigl\{\qb_t^{\ov x},\quad
t\in\ZZ^\nu\text{ and }\ov x\in\cm X^{t\cmp}\bigr\}$ is (the unique)
quasilocal version of one-point conditional distribution of a strictly
positive random field $\pb$. For any $t,s\in\ZZ^\nu$, $x,u\in\cm X^t$ and
$y,v\in\cm X^s$ consider the function $\ov x\mapsto\qb_t^{\ov x
v}(x)\qb_{s\zt}^{\ov x x}(y)\qb_t^{\ov x y}(u)\qb_{s\zt}^{\ov x
u}(v)-\qb_{s\zt}^{\ov x u}(y)\qb_t^{\ov x y}(x)\qb_{s\zt}^{\ov x
x}(v)\qb_t^{\ov x v}(u)$ on $\cm X^{\{t,s\}\cmp}$. Applying the reasoning used
in the proof of Proposition~\ref{P1}, it clearly comes that this function is
equal to zero everywhere.
\end{proof}

Concluding this section, let us note that combining the definition of Gibbs
random field with Criterion~\ref{C3} and taking into account the results of
Section~\ref{OneSpec}, one gets the following characterization.

\begin{criterion}
\label{C6}
A random field is a Gibbs random field if and only if it has a version of
one-point conditional distribution which can be expressed by Gibbs
formulae~\eqref{gibbs} and~\eqref{energy} with some uniformly convergent
potential.
\end{criterion}

\subsection{Random field Gibbsianness criteria\\
            in terms of one-point finite-conditional distribution~I}
\label{MainOut1}

Now we can establish random field Gibbsianness criteria in terms of one-point
finite-conditional distribution, which are precisely the main outcome of the
present paper. The first such criterion is the following.

\begin{criterion}
\label{C7}
A random field is a Gibbs random field if and only if it is strictly positive
and its one-point finite-conditional distribution\/ $\bigl\{\bm q_t^{\fov
x},\quad t\in\ZZ^\nu\text{ and\/ }\fov x\in\wt{\cm X^t}\bigr\}$ satisfy one of
the following equivalent conditions:
\begin{itemize}
\item[\bf (A)]
the limits
$$
\Mlim_{\Lambda\uparrow\ZZ^\nu\sm t}\bm q_t^{\ov x_{\Lambda}}(x)\ ,\quad
t\in\ZZ^\nu,\ x\in\cm X^t,\ \ov x\in\cm X^{t\cmp},
$$
exist, are nonnull uniformly with respect to $\ov x$, and the convergence is
uniform with respect to\/~$\ov x$,
\item[\bf (B)]
the limits
$$
\Mlim_{\Lambda\uparrow\ZZ^\nu\sm t}\bm q_t^{\ov x_{\Lambda}}(x)\ ,\quad
t\in\ZZ^\nu,\ x\in\cm X^t,\ \ov x\in\cm X^{t\cmp},
$$
exist, are strictly positive, and the convergence is uniform with respect to\/
$\ov x$.
\end{itemize}
\end{criterion}

\begin{proof}
The sufficiency is quite evident. Indeed, the strictly positive limits
supposed to exist form a strictly positive version of one-point conditional
distribution of the random field. The uniformity of convergence guarantees
that this version is quasilocal and so, the sufficiency follows from
Criterion~\ref{C5}. Let us also note that at the same time this quasilocality
clearly yields the equivalence of the conditions~(A) and~(B).

Now let us turn to the proof of the necessity. Let $\pb$ be a Gibbs random
field. According to Criterion~\ref{C5} it has a quasilocal and strictly
positive version $\bca Q=\bigl\{\qb_t^{\ov x},\quad t\in\ZZ^\nu\text{ and\/
}\ov x\in\cm X^{t\cmp}\bigr\}$ of one-point conditional distribution. So, to
conclude the proof it is sufficient to show that
$$
\Mlim_{\Lambda\uparrow\vphantom{\cm X^{t\cmp}}\ZZ^\nu\sm t}\ \ \sup_{\ov
x\in\cm X^{t\cmp}\vphantom{\cm X^{t\cmp}}}\ \ \left|\bm q_t^{\ov
x_\Lambda}(x)-\qb_t^{\ov x}(x)\right|=0
$$
for all $t\in\ZZ^\nu$ and $x\in\cm X^t$.

For this we need the following inequality due to Sullivan:
\begin{equation}
\label{ineqsull}
\inf_{\ov y\in\cm X^{t\cmp}\st\ov y_\Lambda=\bm z}\qb_t^{\ov y}(x)
\leq\bm q_t^{\bm z}(x)\leq 
\sup_{\ov y\in\cm X^{t\cmp}\st\ov y_\Lambda=\bm z}\qb_t^{\ov y}(x)
\end{equation}
for all $t\in\ZZ^\nu$, $\Lambda\in\cm E^*(t\cmp)$, $x\in\cm X^t$ and $\bm
z\in\cm X^\Lambda$. This inequality is clearly valid since
$$
\bm q_t^{\bm z}(x)= \frac{\pb_{t\cup\Lambda}(x\bm z)}{\pb_\Lambda(\bm z)}=
\frac1{\pb_\Lambda(\bm z)} \int\limits_{\{\ov y\in\cm X^{t\cmp}\st\ov
y_\Lambda=\bm z\}}\qb_t^{\ov y}(x)\;\pb_{t\cmp}(d\ov y).
$$

Taking this inequality into account, it remains to verify that
$$
\Mlim_{\Lambda\uparrow\vphantom{\cm X^{t\cmp}}\ZZ^\nu\sm t}\ \ \sup_{\ov
x\in\cm X^{t\cmp}\vphantom{\cm X^{t\cmp}}}\ \ \left| \inf_{\ov y\in\cm
X^{t\cmp}\st\ov y_\Lambda=\ov {x\vphantom{y}}_\Lambda}\qb_t^{\ov y}(x)
-\qb_t^{\ov x}(x)\right|=0
$$
and
$$
\Mlim_{\Lambda\uparrow\vphantom{\cm X^{t\cmp}}\ZZ^\nu\sm t}\ \ \sup_{\ov
x\in\cm X^{t\cmp}\vphantom{\cm X^{t\cmp}}}\ \ \left| \sup_{\ov y\in\cm
X^{t\cmp}\st\ov y_\Lambda=\ov {x\vphantom{y}}_\Lambda}\qb_t^{\ov y}(x)
-\qb_t^{\ov x}(x)\right|=0
$$
for all $t\in\ZZ^\nu$ and $x\in\cm X^t$. To show the first one we write
$$
\sup_{\ov x\in\cm X^{t\cmp}}\ \ \left|\inf_{\ov y\in\cm X^{t\cmp}\st\ov
y_\Lambda=\ov {x\vphantom{y}}_\Lambda}\qb_t^{\ov y}(x) -\qb_t^{\ov
x}(x)\right|\leq \sup_{\ov x\in\cm X^{t\cmp}}\ \ \sup_{\ov y\in\cm
X^{t\cmp}\st\ov y_\Lambda=\ov {x\vphantom{y}}_\Lambda}\ \ \left|\qb_t^{\ov
y}(x) -\qb_t^{\ov x}(x)\right|
$$
and use the quasilocality of $\bca Q$. The second one is proved similarly.
\end{proof}

Roughly speaking, Criterion~\ref{C7} asserts that aside from positivity
consideration, Gibbs random fields are characterized by the uniform
convergence of their one-point finite-conditional distribution (to the
one-point conditional one), while only a weaker (almost sure) convergence is
guaranteed for a general random field. In our opinion, this is perhaps the
most comprehensible characterization of Gibbs random fields, on the basis of
which the following purely probabilistic definition of Gibbs random field can
be given.

\begin{definition}
A random filed\/ $\pb$ is called\/ \emph{Gibbs random field\/} if
\begin{itemize}
\item[1)]
for any\/ $\Lambda\in\cm E$ and\/ $\bm x\in\cm X^\Lambda$ one has\/
$\pb_\Lambda(\bm x)>0$,
\item[2)]
the limits
$$
\Mlim_{\Lambda\uparrow\ZZ^\nu\sm t}\frac{\pb_{t\cup\Lambda}(x\ov
x_{\Lambda})}{\pb_{\Lambda}(\ov x_{\Lambda})}\ ,\quad t\in\ZZ^\nu,\ x\in\cm
X^t,\ \ov x\in\cm X^{t\cmp},
$$
exist, are strictly positive, and the convergence is uniform with respect to\/
$\ov x$.
\end{itemize}
\end{definition}

Taking this definition as the definition of Gibbs random field, the usual one
turns into the following representation theorem.

\begin{T}
If\/ $\pb$ is a Gibbs random field, then\/ $\pb$ has a version of conditional
distribution which can be expressed by Gibbs formulae~\eqref{gibbs}
and~\eqref{energy} with some uniformly convergent potential.

Conversely, if a random field\/ $\pb$ has a version of conditional
distribution which can be expressed by Gibbs formulae~\eqref{gibbs}
and~\eqref{energy} with some uniformly convergent potential, then\/ $\pb$ is a
Gibbs random field.
\end{T}

\subsection{Random field Gibbsianness criteria\\
            in terms of one-point finite-conditional distribution~II}
\label{MainOut2}

At first sight, the above presented Criterion~\ref{C7} deals only with
one-point finite-conditional distribution. However, in fact it imposes
conditions equally on its limit, that is, on one-point conditional
distribution. The following and last criterion really deals only with
one-point finite-conditional distribution. Before formulating it, let as agree
that in the sequel when we use the notation $\bm x_T$ we presume that only
configurations~$\bm x$ such that $\supp(\bm x)\supset T$ are considered.

\begin{criterion}
\label{C8}
A random field is a Gibbs random field if and only if it is strictly positive,
its one-point finite-conditional distribution\/ $\bigl\{\bm q_t^{\fov x},\quad
t\in\ZZ^\nu\ \text{and\/}\ \fov x\in\wt{\cm X^t}\bigr\}$ is uniformly nonnull
and one of the following equivalent conditions holds:
\begin{itemize}
\item[\bf (C)]
for any\/ $t\in\ZZ^\nu$ and\/ $x\in\cm X^t$ one has
$$
\Mlim_{\Lambda\uparrow\vphantom{\wt{\cm X^t}}\ZZ^\nu\sm t}\ \ \sup_{\fov
x,\fov y\in\wt{\cm X^t}\st\fov {x\vphantom{y}}_\Lambda=\fov y_\Lambda}\ \
\left|\bm q_t^{\fov x\vphantom{\fov y}}(x)-\bm q_t^{\fov y}(x)\right|=0,
$$
\item[\bf (D)]
for any\/ $t\in\ZZ^\nu$ and\/ $x\in\cm X^t$ one has
$$
\Mlim_{\Lambda\uparrow\vphantom{\wt{\cm X^t}}\ZZ^\nu\sm t}\ \ \sup_{J\in
\vphantom{\wt{\cm X^t}}\cm E^*(t\cmp)}\ \ \sup_{\fov x,\fov y\in\cm X^J
\vphantom{\wt{\cm X^t}}\st\fov {x\vphantom{y}}_\Lambda=\fov y_\Lambda}\ \
\left|\bm q_t^{\fov x\vphantom{\fov y}}(x)-\bm q_t^{\fov y}(x)\right|=0,
$$
\item[\bf (E)]
for any\/ $t\in\ZZ^\nu$ and\/ $x\in\cm X^t$ one has
$$
\Mlim_{\Lambda\uparrow\vphantom{\wt{\cm X^t}}\ZZ^\nu\sm t}\ \ \sup_{\fov
x\in\wt{\cm X^t}}\ \ \left|\bm q_t^{\fov x\vphantom{\fov x_\Lambda}}(x)-\bm
q_t^{\fov x_\Lambda}(x)\right|=0.
$$
\end{itemize}
\end{criterion}

\begin{proof}
First we concentrate on the condition~(E). Clearly
$$
\sup_{\fov x\in\wt{\cm X^t}}\ \ \left|\bm q_t^{\fov x\vphantom{\fov
x_\Lambda}}(x)-\bm q_t^{\fov x_\Lambda}(x)\right|=\sup_{I\in\vphantom{\wt{\cm
X^t}}\cm E(t\cmp)\st I\supset\Lambda} \ \ \sup_{\ov x\in\vphantom{\wt{\cm
X^t}}\cm X^{t\cmp}}\ \ \left|\bm q_t^{\ov x_I}(x)-\bm q_t^{\ov
x_\Lambda}(x)\right|,
$$
and so the condition~(E) is nothing but the Cauchy condition for the existence
of the uniform limits considered in Criterion~\ref{C7}. The sufficiency now
clearly follows from Criterion~\ref{C7} since the Cauchy principle yields the
existence of the uniform limits, and the uniform nonnullness of one-point
finite-conditional distribution guarantees their strict positivity. The
necessity also follows from Criterion~\ref{C7} since the condition~(E) is
ensured by the Cauchy principle, and the uniform nonnullness of one-point
finite-conditional distribution can be easily obtained from~\eqref{ineqsull}
and the condition~(A) $\bigl($use the first inequality of~\eqref{ineqsull} and
the uniform nonnullness of limits considered in the condition~(A)$\bigr)$.

It remains to check the equivalence of the conditions~(C),~(D) and~(E). The
implications (C)$\Rightarrow$(D) and (C)$\Rightarrow$(E) are trivial since
$$
\sup_{J\in \vphantom{\wt{\cm X^t}}\cm E^*(t\cmp)}\ \ \sup_{\fov x,\fov y\in\cm
X^J \vphantom{\wt{\cm X^t}}\st\fov {x\vphantom{y}}_\Lambda=\fov y_\Lambda}\ \
\left|\bm q_t^{\fov x\vphantom{\fov y}}(x)-\bm q_t^{\fov y}(x)\right|
\leq\sup_{\fov x,\fov y\in\wt{\cm X^t}\st\fov {x\vphantom{y}}_\Lambda=\fov
y_\Lambda}\ \ \left|\bm q_t^{\fov x\vphantom{\fov y}}(x)-\bm q_t^{\fov
y}(x)\right|
$$
and
$$
\sup_{\fov x\in\wt{\cm X^t}}\ \ \left|\bm q_t^{\fov x\vphantom{\fov
x_\Lambda}}(x)-\bm q_t^{\fov x_\Lambda}(x)\right|\leq \sup_{\fov x,\fov
y\in\wt{\cm X^t}\st\fov {x\vphantom{y}}_\Lambda=\fov y_\Lambda}\ \ \left|\bm
q_t^{\fov x\vphantom{\fov y}}(x)-\bm q_t^{\fov y}(x)\right|.
$$
Similarly, the inequality
$$
\sup_{\fov x,\fov y\in\wt{\cm X^t}\st\fov {x\vphantom{y}}_\Lambda=\fov
y_\Lambda}\ \ \left|\bm q_t^{\fov x\vphantom{\fov y}}(x)-\bm q_t^{\fov
y}(x)\right| \leq 2\sup_{\fov x\in\wt{\cm X^t}}\ \ \left|\bm q_t^{\fov
x\vphantom{\fov x_\Lambda}}(x)-\bm q_t^{\fov x_\Lambda}(x)\right|
$$
yields the implications (E)$\Rightarrow$(C). To prove the last implication
(D)$\Rightarrow$(C), we need the following lemma.

\begin{lemma}
Let\/ $\bigl\{\bm q_I^{\fov x},\quad I\in\cm E\ \text{and\/}\ \fov x\in\wt{\cm
X^I}\bigr\}$ be the finite-conditional distribution of some strictly positive
random field. Then the set
$$
A=\Bigl\{\ov x\in\cm X^{\ZZ^\nu}\st\Mlim_{\Lambda\uparrow\ZZ^\nu\sm I}\bm
q_I^{\ov x_{\Lambda}}(x)\ \text{exists for every\/}\ I\in\cm E\ \text{and\/}\
x\in\cm X^I\Bigr\}
$$
is of probability\/~$1$ and possesses the following property: if\/ $\ov x\in
A$ then\/ $\bm z\ov x_{J\cmp}\in A$ for all\/ $J\in\cm E$ and\/ $\bm z\in\cm
X^J$.
\end{lemma}

\begin{proof}
Since the set $A$ is a countable intersection of sets of probability~$1$, it
is also of probability~$1$. It remains to show that if $\ov x\in A$ then $\ov
y=z\ov x_{t\cmp}\in A$ for all $t\in\ZZ^\nu$ and $z\in\cm X^t$, that is,
$\Mlim\limits_{\Lambda\uparrow\ZZ^\nu\sm I}\bm q_I^{\ov y_{\Lambda}}(\bm x)$
exists for every $I\in\cm E$ and $\bm x\in\cm X^I$. This is trivial if $t\in
I$ (since in this case $\ov y_{\Lambda}=\ov x_{\Lambda}$) and clearly follows
from the relation
$$
\bm q_I^{\ov y_{\Lambda}}(\bm x)=\bm q_I^{z\ov x_{\Lambda\sm t}}(\bm x)
=\frac{\bm q_{t\cup I}^{\ov x_{\Lambda\sm t}}(z\bm x)} {\bigl(\bm q_{t\cup
I}^{\ov x_{\Lambda\sm t}}\bigr)_t(z)},\quad\Lambda\ni t,
$$
otherwise.
\renewcommand{\qedsymbol}{$\vartriangleleft$}
\end{proof}

Returning to the proof of the implication (D)$\Rightarrow$(C), let us fix some
$t\in\ZZ^\nu$ and $x\in\cm X^t$, denote
$$
f(\Lambda)=\sup_{J\in \vphantom{\wt{\cm X^t}}\cm E^*(t\cmp)}\ \ \sup_{\fov
x,\fov y\in\cm X^J \vphantom{\wt{\cm X^t}}\st\fov {x\vphantom{y}}_\Lambda=\fov
y_\Lambda}\ \ \left|\bm q_t^{\fov x\vphantom{\fov y}}(x)-\bm q_t^{\fov
y}(x)\right|,
$$
and for any $\varepsilon>0$ choose $\bigl($according to the condition
(D)$\bigr)$ some $\Lambda_\varepsilon\in\cm E$ such that
$\left|f(\Lambda)\right|<\varepsilon$ for all $\Lambda\in\cm E$,
$\Lambda\supset\Lambda_\varepsilon$.

First, we will show that $\Mlim\limits_{\Lambda\uparrow\ZZ^\nu\sm t}\bm
q_t^{\ov x_{\Lambda}}(x)$ exists for every $\ov x\in\cm X^{t\cmp}$. Let us
take some $\ov x^\circ\in A$ (according to the lemma, the set $A$ is of
probability~1 and so is not empty) and consider $\ov y=\ov
x^{\vphantom{\circ}}_{\Lambda_\varepsilon}\ov
x^\circ_{\Lambda_\varepsilon\cmp}\in A$. So, we can find some
$\Lambda_\varepsilon'\in\cm E$ such that $\left|\bm q_t^{\ov y_I}(x)-\bm
q_t^{\ov y_J}(x)\right|<\varepsilon$ for all $I,J\in\cm E$,
$I\supset\Lambda_\varepsilon'$, $J\supset\Lambda_\varepsilon'$. Thus, for all
$I,J\in\cm E$ such that $I\supset\Lambda_\varepsilon\cup\Lambda_\varepsilon'$
and $J\supset\Lambda_\varepsilon\cup\Lambda_\varepsilon'$ we can write
\begin{align*}
\left|\bm q_t^{\ov x_I}(x)-\bm q_t^{\ov x_J}(x)\right|
&\leq\left|\bm q_t^{\ov x_I}(x)-\bm q_t^{\ov y_I}(x)\right|
+\left|\bm q_t^{\ov y_I}(x)-\bm q_t^{\ov y_J}(x)\right|
+\left|\bm q_t^{\ov y_J}(x)-\bm q_t^{\ov x_J}(x)\right|\\
&<f(\Lambda_\varepsilon)+\varepsilon+f(\Lambda_\varepsilon)<3\varepsilon,
\end{align*}
and hence $\Mlim\limits_{\Lambda\uparrow\ZZ^\nu\sm t}\bm q_t^{\ov
x_{\Lambda}}(x)$ exists according to Cauchy principle.

Further, for every $\ov x\in\cm X^{t\cmp}$ consider the set $V(\ov
x)=\left\{\ov y\in\cm X^{t\cmp}\st\ov y_{\Lambda_\varepsilon}=\ov
{x\vphantom{y}}_{\Lambda_\varepsilon}\right\}$. Clearly these sets are either
mutually disjoint or coinciding, and there is only a finite number~$k$
$\bigl($more precisely $k=\bil|\cm X^{\Lambda_\varepsilon}\bir|\bigr)$ of
different sets among them.  Hence there exists a finite collection $\ov
x^1,\ldots,\ov x^k\in\cm X^{t\cmp}$ such that
$$
\cm X^{t\cmp}=\bigcup_{i=1}^{k}V(\ov x^i).
$$
$\bigl($This fact equally follows from the compactness of $\cm
X^{t\cmp}$.$\bigl)$ So, using Cauchy principle we can find some
$\Lambda_\varepsilon''\in\cm E$ such that $\left|\bm q_t^{\ov x^i_I}(x)-\bm
q_t^{\ov x^i_J}(x)\right|<\varepsilon$ for all $i=1,\ldots,k$ and all
$I,J\in\cm E$, $I\supset\Lambda_\varepsilon''$,
$J\supset\Lambda_\varepsilon''$.

Now, let the set $\Lambda\in\cm E$ be such that
$\Lambda\supset\Lambda_\varepsilon\cup\Lambda_\varepsilon''$, the sets
$I,J\in\cm E$ be such that $I\supset\Lambda$ and $J\supset\Lambda$, and the
configurations $\fov x\in\cm X^I$ and $\fov y\in\cm X^J$ be such that $\fov
{x\vphantom{y}}_\Lambda=\fov y_\Lambda$. Clearly, we can find some
$i\in\{1,\ldots,k\}$ such that $\ov x^i_{\Lambda_\varepsilon}=\fov
{x\vphantom{y}}_{\Lambda_\varepsilon}=\fov y_{\Lambda_\varepsilon}$, and thus
we may write
\begin{align*}
\left|\bm q_t^{\fov x\vphantom{\fov y}}(x)-\bm q_t^{\fov y}(x)\right|&\leq
\left|\bm q_t^{\fov x\vphantom{\fov y}}(x)-\bm q_t^{\ov x^i_I}(x)\right|
+\left|\bm q_t^{\ov x^i_I}(x)-\bm q_t^{\ov x^i_J}(x)\right|
+\left|\bm q_t^{\ov x^i_J}(x)-\bm q_t^{\fov y}(x)\right|\\
&<f(\Lambda_\varepsilon)+\varepsilon+f(\Lambda_\varepsilon)<3\varepsilon
\end{align*}
which shows that the condition (C) holds.
\end{proof}

In conclusion let us note that the analogues of Criteria~\ref{C7} and~\ref{C8}
formulated in terms of all the finite-conditional distribution are of course
valid. Concerning the first one, we would like to mention that its necessity
statement was as a matter of fact contained in the proof of Lemma~1 of
Sullivan~\cite{sull2}, whose argument we follow while proving
Criterion~\ref{C7}. As to the second one, let us mention that the part
utilizing the analogue of the condition~(D) can be deduced from Theorems~1
and~2 of Kozlov~\cite{kozl}. It should be pointed out that the author does not
provide the proof of the sufficiency statement of Theorem~2 (leaving it, as he
says, to the reader). However, our considerations show that the proof of this
statement is neither intuitive, nor technically simple. Moreover, the validity
of the statement seems dubitable in the settings of Kozlov~\cite{kozl} where
the state space is not supposed to be finite or even compact.


\end{document}